
\documentstyle{amsppt}
\magnification=1200
\hyphenation{}
\define\Z{{\Bbb Z}}

\define\Q{{\Bbb Q}}
\define\M{\Cal{M}}
\define\J{\Cal{J}}
\define\G{{\Gamma}}
\define\C{{\Bbb C}}
\define\ra{\rangle}
\define\la{\langle}
\define\Cu{\Cal{C}}
\define\Sy{{\frak S}}
\define\Hom{\operatorname{Hom}}

\define\Sym{\operatorname{Sym}}
\define\Pic{\operatorname{Pic}}

\define\symp{\operatorname{Sp}}

\define\codim{\operatorname{codim}}

\define\refer#1{(#1)}

\newcount\headnumber
\newcount\labelnumber
\define\section{\global\advance\headnumber
by1\global\labelnumber=0{{\the\headnumber}.\ }}
\define\label{(\global\advance\labelnumber by1 \the\headnumber
.\the\labelnumber )\enspace}

\NoBlackBoxes

\topmatter

\title
Stable cohomology of the mapping class group with symplectic coefficients
and of the universal Abel--Jacobi map
\endtitle

\rightheadtext{Stable cohomology of mapping class groups}
\author
Eduard Looijenga
\endauthor

\abstract
The irreducible representations of the complex symplectic group of genus $g$
are indexed by nonincreasing sequences of integers $\lambda =(\lambda _1\ge
\lambda_2\ge\cdots)$ with $\lambda _k=0$ for $k>g$. A
recent result of N.V.~Ivanov implies that for a given partition $\lambda$, the
cohomology group of a given degree of the mapping class group of genus $g$
with values in the representation associated to $\lambda$ is independent of $g$
if $g$ is sufficiently large. We prove that this stable cohomology is the
tensor product of the stable cohomology of the mapping class group and a
finitely generated graded module over $\Q [c_1,\dots ,c_{|\lambda |}]$, where
$\deg (c_i)=2i$ and $|\lambda |=\sum _i\lambda _i$. We describe this module
explicitly. In the same sense we determine the
stable rational cohomology of the moduli space of compact Riemann surfaces with
$s$ given ordered distinct (resp. not necessarily distinct) points as well as
the stable cohomology of the universal Abel--Jacobi map. These results take
into account mixed Hodge structures.
\endabstract

\address
Faculteit Wiskunde en Informatica,
Universiteit Utrecht,
PO Box 80.010, 3508 TA Utrecht,
The Netherlands
\endaddress

\email
looijenga\@math.ruu.nl
\endemail



\endtopmatter

\document

\head
\section Introduction
\endhead

The mapping class group $\G _{g,r}^s$ can be defined in terms of a compact
connected oriented surface $S_g$ of genus $g$ on which are given $s+r$
(numbered) distinct points $(x_i)_{i=1}^{r+s}$: it is then the connected
component group of the group of orientation preserving diffeomorphisms of $S_g$
which fix each $x_i$ and are the identity on the tangent space of $S_g$ at
$x_i$ for $i=s+1,\dots ,s+r$. It is customary to omit the suffix $r$ resp. $s$
when it is zero.

Harer's stability theorem says essentially that $H^k(\G _{g,r}^s;\Z )$ only
depends on $s$ if $g$ is large compared to $k$.  For a more precise statement
it
is convenient to make a definition first. There is a natural outer homomorphism
$\Gamma ^s_{g,r+1}\to \Gamma ^s_{g+1,r}$ (that is, an orbit of homomorphisms
under the  inner automorphism group of the target group, so that there is
well-defined map on homology) and there is a forgetful homomorphism $\Gamma
^s_{g,r+1}\to \Gamma ^s_{g,r}$. For a coefficient
ring $R$ and an integer $g_0\ge 0$, we define $N(g_0;R)$ as the maximal integer
$N$ such that both induce isomorphisms on homology with coefficients in $R$ in
degree $\le N$ for all $g\ge g_0$ and $s,r\ge 0$. Harer showed in \cite{3} that
$N(g;\Z )\ge {1\over 3}g$ and Ivanov \cite{5}, \cite{6} improved this to
$N(g;\Z
)\ge {1\over 2}g -1$. Recently, Harer proved that $N(g;\Q )\le {2\over 3}g$
and that $N(g;\Q )\ge {2\over 3}g$ is almost true: it holds, provided that we
restrict to the mapping class groups with $r\ge 1$ \cite{4}. It is likely that
in fact
$N(g;\Q )\ge {2\over 3}g-1$. We will be mostly concerned with $N(g;\Q )$ and so
we shall write for this number $N(g)$ instead. A consequence of the stability
property is that for every integer $s\ge 0$ we have a stable cohomology algebra
$H^{\bullet}(\G _{\infty}^s;R)$ (where as before, we omit the superscript $s$
if
it is equal to $0$). As the notation suggests, this is indeed
the cohomology of a group $\G _{\infty}^s$, namely the group of compactly
supported mapping classes of an oriented connected surface of infinite genus
relative $s$ given numbered points. (The number of ends of this surface may be
arbitrary.)

\medskip
Consider the symplectic vector
space $V_g:=H^1(S_g;\Q )$. Its symplectic form is preserved by the natural
action of the mapping class group $\G _g$ on $V_g$ and so any finite
dimensional
representation $U$ of the algebraic group  $Sp(V_g)$ can be regarded as a
representation of $\G _g$; in particular we have defined the cohomology groups
$H^k(\G _g ;U)$. A basic fact of representation theory is that the isomorphism
classes of the irreducible complex representations of $Sp(V_g)$ are in a
natural bijective correspondence with $g$-tuples of nonnegative integers
$(a_1,\dots ,a_g)$. For instance, the $k$th basis vector $(0,\dots ,0,1,0\dots
,0)$ corresponds to the $k$th exterior power of $V_g$. This result goes back to
Weyl, who gave in addition a functorial construction of such a
representation inside  $\Sym ^{a_1}(\wedge ^1 V_g)\otimes\cdots\otimes\Sym
^{a_g}(\wedge ^g V_g)$. He labeled this representation by the sequence
$(a_1+\cdots +a_g,a_2+a_3+\cdots +a_g,\dots ,a_{g-1}+a_g,a_g)$. We will follow
his convention, so for us a numerical partition $\lambda=(\lambda _1\ge \lambda
_2\ge\cdots )$ with at most $g$ parts (that is, $\lambda _k=0$ for $k>g$)
determines an
irreducible representation  $S_{\langle {\lambda}\rangle}(V_g)$. It follows
from
a recent result of Ivanov \cite{6} that for fixed $k$ and $\lambda$, the
cohomology groups $H^k(\G _{g,1};S_{\langle {\lambda}\rangle}(V_g))$ are
independent of $g$ if $k\le N(g)-|\lambda |$, where $|\lambda |:=\sum _i\lambda
_i$ is the size of the partition.  We shall give an
independent proof of this in the undecorated case $r=s=0$ which
generalizes in a straightforward manner to the case of arbitrary $r$ and $s$
and we determine these stable cohomology groups as $H^{\bullet}(\G
_{\infty} ;\Q)$-modules at the same time.

\proclaim
{\label Theorem} For every numerical partition $\lambda$ of $s$, there is
a graded finitely generated $\Q [c_1,\dots
,c_s]$-module $B^{\bullet}_{\lambda}$ (where $c_i$ has degree $2i$) and a
natural homomorphism
$$
H^{\bullet}(\G _{\infty};\Q )\otimes B^{\bullet}_{\lambda}\to
H^{\bullet}(\G _g; S_{\langle{\lambda}\rangle}(V_g)),
$$
which is an isomorphism in degree $\le N(g)-|\lambda |$.
\endproclaim

To describe $B^{\bullet}_{\lambda}$, denote the coordinates of $s$-space
$\Q ^s$ by $u_1,\dots ,u_s$ so that $\Q [u_1,\dots ,u_s]$ is its algebra of
regular functions. We grade this algebra by giving each coordinate $u_i$ weight
$2$. A {\sl diagonal} of $\Q ^s$ is by definition  an intersection of the
hyperplanes $u_i=u_j$; this includes the intersection with empty index set,
that
is, $\Q ^s$ itself. It is clear that a partition $P$ of the set $\{ 1,\dots
,s\}$
determines (and is determined by) a diagonal $\Delta _P$. Notice that the
algebra of regular functions $\Q [\Delta _P]$ is a quotient of
$\Q [u_1,\dots ,u_s]$ by a graded ideal, so that it inherits a grading.

Denote by $l_i(P)$ the number of parts of $P$ of cardinality $i$, and
consider
$$
\bigoplus _P t^{-s}u^{\codim (P)+2l_1(P)+l_2(P)}\Q [\Delta _P],\tag{1}
$$
where $\codim (P)$ is short for $\codim (\Delta _P)$ and $t$ resp. $u$ has
formal degree $1$ resp. $2$. (The difference between
$t^2$ and $u$ becomes manifest in the Hodge theory: $t$ has
Hodge type $(0,0)$, whereas $u$ has Hodge type $(1,1)$.) We regard this as a
graded
module over $\Q [u_1,\dots ,u_s]$ that has a natural action of the symmetric
group $\Sy _s$. We tensorize this module with the signum representation of $\Sy
_s$ and denote the resulting graded  $\Q [u_1,\dots ,u_s]$-module with $\Sy
_s$-action by $B^{\bullet}_s$. Now recall that every numerical partition
$\lambda$ of $s$
determines an equivalence class $(\lambda )$ of irreducible representations
of $\Sy _s$ and that this gives a bijection between these two sets. For
instance, the coarsest partition $(s)$ labels the trivial representation,
whereas the finest partition $(1^s)$ corresponds to the signum representation.
Passing from a partition $\lambda$ to the conjugate partition $\lambda '$
corresponds to taking the tensor product with the signum representation. In the
case at hand we have a decomposition of $B^{\bullet}_s$ into isotypical
subspaces:
$$
B^{\bullet}_s=\oplus _{\lambda} B^{\bullet}_{\lambda}\otimes (\lambda),\quad
\text{with } B^{\bullet}_{\lambda}=\Hom _{\Sy _s}((\lambda ),B^{\bullet}_s).
$$
Clearly, this is also a decomposition into graded $\Q
[u_1,\dots ,u_s]^{\Sy _s}$-submodules. We identify the latter ring with
$\Q [c_1,\dots ,c_s]$, where $c_k$ is the $k$th elementary symmetric function
in the $u_i$'s. This completes the description of $B^{\bullet}_{\lambda}$.

It follows from the work of M. Saito that $H^{\bullet}(\G _g;
S_{\langle{\lambda}\rangle}(V_g))$ carries a natural mixed Hodge structure.  It
is known that $H^{\bullet}(\G _{\infty};\Q )$ has a natural mixed Hodge
structure as well (see \refer{2.5}). We will find that if we put a Hodge
structure on $B^{\bullet}_{\lambda}$ by giving its degree $2d-|\lambda |$-part
Hodge type $(d,d)$, then the homomorphism of \refer{1.1} is a morphism of mixed
Hodge structures. M.\ Pikaart has recently shown that $H^n(\G _{\infty};\Q )$
is pure of weight $n$; this implies that $H^n(\G _g; S_{\la\lambda\ra
}(V_g))$ is pure of weight $n+|\lambda |$ in the stable range $n\le
N(g)-|\lambda |$.
\medskip

{\it Example 1.} The $s$-th symmetric power of $V_g$ corresponds to
$S_{\langle s\rangle}(V_g)$ and so the stable cohomology of $\Gamma _g$ with
values in $\Sym ^s(V_g)$ is the tensor product of $H^{\bullet}(\G _{\infty} )$
with
$B^{\bullet}(1^s)$, that is, the isotypical subspace of expression \refer{1}
corresponding to the signum representation of $\Sy _s$. Any partition different
from the partition into singletons is invariant under a transposition and so
will not contribute. This leaves us therefore with the $(1^s)$-isotypical
subspace of $t^{-s}u^{2s}\Q [u_1,\dots ,u_s]$. This is the free $\Q [c_1,\dots
,c_s]$-module generated by the element $t^{-s}u^{2s}\prod _{i>j} (u_i-u_j)$,
hence is of the form $t^{-s}u^{{1\over 2}s(s+3)}\Q [c_1,\dots ,c_s]$. We find
that $$ H^{\bullet}(\G _{\infty };\Sym ^s (V_{\infty}))= t^{-s}u^{{1\over
2}s(s+3)}H^{\bullet}(\G _{\infty };\Q )[c_1,\dots ,c_s]. $$

\smallskip
{\it Example 2.} For $g\ge s$, the primitive subspace $\Pr ^s(V_g)$ of the
$s$-th exterior power of $V_g$ corresponds to
$S_{\langle 1^s\rangle V_g}$ and so the stable cohomology of $\Gamma _g$ with
values in $\Pr ^s(V_g)$ is the tensor product of $H^{\bullet}(\G _{\infty};\Q
)$ with ${\Sy _s}$-invariant part of \refer{1}. This is naturally written as a
sum over the numerical partitions of $s$. Here a
numerical partition is best described by means of the exponential
notation: $(1^{l_1}2^{l_2}3^{l_3}\cdots )$, where $l_k$ is the number of parts
of cardinality $k$ (so that $\sum _k kl_k =s$). Its contribution is then
$$
\otimes _{k\ge 1}t^{-l_k k}u^{l_k\max (2, k-1)}\Q [c_1,c_2,\dots ,c_{l_k}].
$$
If we sum over all sequences $(l_1,l_2,\dots )$ of nonnegative integers which
become eventually zero, we get
$$
\otimes _{k\ge 1}\bigl( \oplus _{l=0}^{\infty} t^{-lk}u^{l\max (2, k-1)}\Q
[c_1,c_2,\dots ,c_l]\bigr) .
$$
and $B^{\bullet}_{(1^s)}$ is the $t^{-s}$-part of this expression.
For instance,
$$\align
B^{\bullet}_{(1^2)}&= t^{-2}(u^4\Q [c_1,c_2]\oplus u^2\Q [c_1]),\\
B^{\bullet}_{(1^3)})&=t^{-3}(u^6\Q [c_1,c_2,c_3]\oplus u^4\Q [c_1]\otimes\Q
[c_1]\oplus u^2\Q [c_1]).
\endalign
$$
Since $H^1(\G _{\infty};\Q )=0$, it follows in particular that $H^1(\G
_{\infty},S_{\langle 1^3\rangle})$ has dimension one.  (It is easy to see that
$H^1(\G _{\infty};,S_{\langle \lambda\rangle})=0$ for all other numerical
partitions $\lambda$.)

\medskip\label
The  cohomology groups of $\G _g$ with values in a symplectic representation
have a geometric interpretation as the cohomology of a local system (in fact,
of
a variation of polarized Hodge structure). The Teichm\"uller space $\Cal{T}_g$
of conformal structures on $S_g$ modulo isotopy is a {\it contractible} complex
manifold (of complex dimension $3g-3$, when $g\ge 2$). The action of $\G _g$ on
it is properly discontinuous and a subgroup of finite index acts freely. The
orbit space $\M _g:=\G _g \backslash \Cal{T}_g$  is naturally interpreted as
the
coarse moduli space of smooth complex projective curves. Via this
interpretation, it gets the structure of a normal quasi-projective variety. It
follows from the preceding that  $\M _g$ has the rational cohomology of $\G
_g$.
More generally, if $U$ is a rational representation of $\G _g^s$, then we have
a
natural isomorphism  $$
H^{\bullet}(\G _g;U)\cong H^{\bullet}(\M _g;\Bbb{U}),
$$
where $\Bbb{U}$ is the sheaf on $\M _g$ which is the
quotient of the trivial local system $\Cal{T}_g\times U\to \Cal{T}_g$ by the
(diagonal) action of $\G _g$. This applies in particular to the
representations $S_{\lambda}(V_g)$. This representation appears in the
cohomology of degree $s:=|\lambda |$ of the configuration space of $s$ numbered
(not necessarily distinct) points on $S_g$. So if $\Cu _g^s\to\M _g$ denotes
the
$s$-fold fiber product of the universal curve, then its $s$th direct image
contains the local system (in the orbifold sense) associated to
$S_{\lambda}(V_g)$ as a direct summand. By a theorem of Deligne, the Leray
spectral sequence of the forgetful map $\Cu _g^s\to M_g$ degenerates at the
$E_2$-term and thus the stable cohomology of this representation is realized
inside $H^{\bullet}(\Cu ^s_g;\Q )$. This fact will be used in an essential
way.

\remark
{Acknowledgements} I thank Dick Hain for discussions and for suggesting to use
Deligne's degeneration theorem. I am also grateful to him for pointing out that
the way \refer{1.1} was stated  in an earlier version could not be correct. His
student S.\ Kabanov obtained related results that should appear soon.

I am also indebted to the referee for invaluable comments, in particular for a
suggestion for shortening the original proof of \refer{2.3}.
\endremark

\head
\section Stable cohomology of $\M _g^s$ and $\Cu _g^s$
\endhead

We first state and prove an immediate consequence of the stability theorems.
For $s\ge 1$, the class of a (Dehn) twist about $x_s$ generates an
infinite cyclic central subgroup of $\G _{g,r+1}^{s-1}$. The quotient group can
be identified with $\G _{g,r}^s$ and so we have a Gysin sequence
$$
\cdots\to  H^{k-2}(\G _{g,r}^s;\Z ){\buildrel
u_s\cup\over\longrightarrow} H^k(\G _{g,r}^s;\Z)\to H^k(\G
_{g,r+1}^{s-1};\Z)\to\cdots,
$$
where $u_s\in H^2(\G _{g,r}^s;\Z )$ is the
first Chern class. Similarly, $x_i$ determines a first Chern
class $u_i\in  H^2(\G _{g,r}^s;\Z )$ for $i=1,\dots ,s$. These classes are
clearly stable and we shall not make any notational distinction between the
$u_i$'s and their stable representatives.

\proclaim{\label Proposition}
The stable cohomology algebra over of the
mapping  class groups of surfaces with $s$ distinct numbered points is a graded
polynomial algebra on the stable cohomology ring of the absolute mapping class
groups. More precisely, there is a natural $\Sy _s$-equivariant graded
ring homomorphism
$$
H^{\bullet}(\G _{\infty};\Z )[u_1,\dots ,u_s]\to H^{\bullet}(\Gamma
^s_{g,r};\Z),
$$
which is an isomorphism in degree $\le N(g;\Z )$.
In particular,  the rational stable cohomology algebra of the
mapping  class groups of surfaces with a $s$ unlabeled points is a graded
polynomial $H^{\bullet}(\G _{\infty};\Q)$-algebra on the elementary symmetric
functions $c_1,\dots ,c_s$ of the $u_i$'s.
\endproclaim
\demo{Proof}
The composite of the forgetful maps
$\G _{g,r}^{s-1}\to \G _{g,r-1}^s$ and $\G _{g,r-1}^s\to \G _{g,r-1}^{s-1}$
induces an isomorphism on $H^k(-;\Z )$ for large $g$. So in this range the
forgetful map
induces a surjection
$H^k(\Gamma _{g,r-1}^s;\Z )\to H^k(\Gamma _{g,r}^{s-1};\Z )$. If we feed
this in the above Gysin sequence, we get short exact sequence
$$
0\to H^{\bullet -2}(\G _{\infty}^s;\Z ){\buildrel u\cup\over\longrightarrow}
H^{\bullet}(\G _{\infty} ^s;\Z )\to H^{\bullet}(\G _{\infty}^{s-1};\Z )\to 0
$$
so that $H^{\bullet}(\G _{\infty}^s;\Z )\cong H^{\bullet}(\G _{\infty}^{s-1};\Z
)[u]$
as algebra's.  The assertion now follows with induction on $s$.
\enddemo

Let us fix a finite set $X$. We denote by
$\Cu ^X_g$  the moduli space of pairs $(C,x )$ where $C$ is a compact Riemann
surface of genus $g$ and $x:X\to C$ is a map. Let $j:\M ^X_g\subset
\Cu ^X_g$ be the open subset defined by the condition that $x$ be injective.
Just as $\M _g$ is a virtual classifying space for $\G _g$, $\M ^X_g$ is one
for
$\G _g^{|X|}$; in particular, $\G _g^{|X|}$ and $\M ^X_g$  have the same
rational cohomology. This enables us to restate \refer{2.1} in more geometric
terms. Let $\Cu _g\to\M _g$ be the universal curve and denote by $\theta$ its
relative tangent sheaf.
For every $i\in X$, the map $(C,x)\mapsto  x(i)$ defines a projection of $\Cu
^X_g$ onto the $\Cu _g$; denote by $\theta _i$ the pull-back of $\theta $ under
this map.
One easily recognizes the first Chern class of $\tau _i|\M ^X_g$ as the first
Chern class defined above. So \refer{2.1} implies:

\proclaim{\label Proposition} The ring homomorphism
$$
H^{\bullet}(\M _g;\Q)[u_i:i\in X]\to H^{\bullet}(\M _g^X;\Q),
\quad u_i\mapsto c_1(\theta _i)|\M ^X_g.
$$
is an isomorphism in degree $\le N(g)$.
\endproclaim

We will use this proposition to prove that the rational cohomology of $\Cu
^X_g$ also stabilizes.

We begin with attaching to $X$ a graded commutative $\Q [u_i :i\in X]$-algebra.
It is convenient to introduce an auxiliary graded
commutative $\Q$-algebra $\tilde  A_X^{\bullet}$ first. The latter is defined
by the following presentation: for each nonempty subset $I$ of $X$, $\tilde
A_X^{\bullet}$ has a generator $u_I$ of degree two (we also write $u_i$ for
$u_{\{ i\}})$, and these are subject to the relations $u_Iu_J=u_iu_{I\cup J}$
if $i\in I\cap J$.
So if $i\in I$, then $u_Iu_I=u_iu_I$. It is then easy to see that the
$\Q [u_i:i\in X]$-submodule generated by $u_I$ is defined by the relations
$(u_i-u_j)u_I=0$ whenever $i,j\in I$. The
monomials $\prod _I u_I^{r_I}$ for which $I$ runs over the members of a
partition of $X$ form an additive basis of $\tilde  A_X^{\bullet}$. (To make
the
indexing effective, let us agree that we only allow $r_I$ to be zero if $I$ is
a singleton.)

We then let $A^{\bullet}_X$ be the  $\Q [u_i:i\in X]$-subalgebra
of  $\tilde A_X^{\bullet}$ generated by the elements
$a_I:=u_I^{|I|-1}$, where $I$ runs over the subsets of $X$ with at least two
elements. These generators obey the relations
$$\align
u_ia_I&:=u_ja_I \text{ if } i,j\in I,\\ a_Ia_J&:= u_i^{|I\cap J| -1}a_{I\cup J}
\text{ if } i\in I\cap J.
\endalign
$$
and it is easy to see we thus obtain a presentation of $A^{\bullet}_X$ as a
graded
commutative $\Q [u_i :i\in X]$-algebra. Notice that as a $\Q [u_i :i\in
X]$-algebra, $A^{\bullet}_X$ is already generated by the $a_I$'s with $|I|=2$.
For every partition $P$ of $X$ we put $a_P:=\prod _{I\in P; |I|\ge 2} a_I$
(with
the convention that $a_P=1$ if $P$ is the partition into singletons). These
elements generate $A^{\bullet}_X$ as a $\Q [u_i :i\in
X]$-module. In fact,
$$
A_X^{\bullet}=\bigoplus _{P|X} \Q [u_I:I\in P]a_P.
$$
We give $\tilde A_X^{\bullet}$ a (rather trivial) Hodge structure: its degree
$2p$-part has Hodge type $(p,p)$. Clearly, $A_X^{\bullet}$ is then a
Hodge substructure.

\smallskip
Given a partition $P$ of $X$, then the pairs $(C,x: X\to
C)$ for which every member of $P$ is contained in  a fiber of
$x$ define a closed  subvariety $i_P:\Cu _g(P)\subset \Cu ^X_g$. Those
for which $P$ is the partition defined by $x$ make up a Zariski-open subvariety
$\M _g(P)\subset \Cu _g(P)$. Notice
that $\Cal{C}_g(P)$ resp. $\M _g(P)$ can be identified with $\Cu ^{X/P}_g$
resp. $\M ^{X/P}_g$ (where $X/P$ stands for quotient of $X$ by the equivalence
relation defined by $P$, or rather the set of parts of $P$), and that this is a
submanifold in the orbifold sense.

For every nonempty $I\subset X$, let $P_I$ be the partition of $X$ whose
parts are $I$ and the singletons in $X-I$. So the corresponding orbifold
$\Cu _g(P_I)$ has codimension $|I|-1$ in $\Cu _g^X$.

\proclaim{\label Theorem} There is an algebra homomorphism
$$
\phi _g^X: H^{\bullet}(\M _g;\Q )\otimes A_X^{\bullet}\to
H^{\bullet}(\Cu ^X_g;\Q )
$$
that extends the natural homomorphism $H^{\bullet}(\M _g;\Q )\to
H^{\bullet}(\Cu ^X_g;\Q )$, sends $1\otimes u_i$ to $c_1(\theta _i)$ and sends
$1\otimes a_I$ to the Poincar\'e dual of the class of  $\Cu _g(P_I)$. This is
an $\Sy _X$-equivariant
algebra homomorphism and is also a morphism of mixed Hodge structures.
Moreover, $\phi _g^X$ is an isomorphism in degree $\le N(g)$.
 \endproclaim
\demo{Proof}
For the first statement we must show that if $I,J\subset X$ have at least two
elements,
$i\in I\cap J$, $j\in I$ and $P$ is a partition of $X$, then
$$
\align
i_{P!}(1)&= \prod _{I'\in P; |I'|\ge 2}i_{P_{I'}!}(1),\\
c_1(\theta _i)i_{P_I!}(1)&=c_1(\theta _j)i_{P_I!}(1),\\
i_{P_I!}(1)i_{P_J!}(1)&=c_1(\theta _i)^{|I\cap J|-1} i_{P_{I\cup J}!}(1).
\endalign
$$
The first identity is geometrically clear and the second follows from the fact
that $\theta _i$ and $\theta _j$ have isomorphic restrictions to $\Cu
_g(P_I)$.
To derive the last identity, we use the
following lemma (the proof of which is left to the reader):
\enddemo

\proclaim{\label Lemma}
Let $U$ and $V$ be closed complex submanifolds of a complex manifold $M$ whose
intersection $W:=U\cap V$ is also a complex manifold. Suppose that
any tangent vector of $M$ which is tangent to both $U$ and $V$ is tangent to
$W$. Then $i_{U!}(1)i_{V!}(1)=i_{W!}(e)$, where $e$ is the
euler class of the cokernel of the natural monomorphism
$\nu _W\to \nu _U|W\oplus\nu _V|W$.
\endproclaim

\demo{Completion of the proof of \refer{2.3}} We apply the orbifold version
of this lemma to $M=\Cu _g^X$, $U=\Cu _g (P_I)$, $V=\Cu _g(P_J)$ so that $W=\Cu
_g(P_{I\cup J})$. The desired assertion then follows if we use the fact that
the
bundles appearing in the monomorphism $\nu _W\to \nu _U|W\oplus\nu _V|W$ are
all
direct sums of copies of $\theta _i|W$.

The second statement of the theorem is clear. To prove the last,
let $U_k$ resp. $S_k$ denote the union of the strata
$\M _g(P)$ of $\codim \le k$ resp.~$=k$. We prove with induction on $k$ that
the homomorphism
$$
\bigoplus _{\codim P\le k}  H^{\bullet}(\M _g;\Q )\otimes \C [u_I:I\in
P]a_P\to  H^{\bullet}(U_k;\Q )
$$
is an isomorphism in degree $\le N(g)$. For $k=0$ this is \refer{2.2}. If $k\ge
1$, then consider the Gysin sequence of the pair $(U_k,S_k)$:
$$
\cdots\to H^{n-2k}(S_k;\Q )\to H^n(U_k;\Q )\to H^n(U_{k-1};\Q )\to\cdots .
$$
In degrees $\le N(g)$ the isomorphism  of $\bigoplus _{\codim P\le k}
H^{\bullet}(\M _g;\Q )\otimes \C [u_I:I\in P]a_P$ onto  $H^n(U_{k-1};\Q )$
factorizes over $H^n(U_k;\Q )$. So the Gysin sequence splits in this range.
Since  $\bigoplus _{\codim P =k}  H^{\bullet}(\M _g;\Q )\otimes \C [u_I:I\in
P]$ maps isomorphically onto $H^{\bullet}(S_k;\Q )$ in degree $\le N(g)$, the
theorem follows.
\enddemo

{\it Remark.}
For a curve $C$, the image of $A^{\bullet}_X$ in $H^{\bullet}(C^X;\Q )$ is
contained in the  Hodge ring of $C^X$. If $C$ is general, then it is in fact
equal to it.

\medskip\label
A virtual classifying space of $\Gamma _g^{r+s}$ is the moduli space of
$r+s$-pointed curves $\M _g^{r+s}$. It carries $r+s$ relative tangent bundles
$\theta _1,\dots ,\theta _{r+s}$ and the total space of the $(\C ^*)^r$-bundle
$\M _{g,r}^s\to \M_g^{r+s}$ defined by  $\theta _{1+s},\dots ,\theta _{r+s}$ is
a virtual classifying space for $\Gamma _{g,r}^s$. The comparison maps that
enter in the statement of the stability theorem  admit simple descriptions in
these algebro-geometric terms and so preserve Hodge structures. This is
probably well-known, but since we do not know a reference, we explain this.
Clearly, the forgetful homomorphism $\Gamma _{g,r+1}^s\to \Gamma _{g,r}^s$
corresponds to the obvious projection $\M _{g,r+1}^s\to \M _{g,r}^s$.  To
exhibit the outer homomorphism
$\Gamma _{g,r+1}^s\to \Gamma _{g+1,r}^s$, we first note that $\M
_g^{r+s+1}\times \M _1^1$ parametrizes a codimension one stratum of the
Knudsen--Deligne--Mumford compactification of $\M _{g+1}^{r+s}$: a pair
$((C;x_1,\dots ,x_{r+s+1}),(E;O))$ determines a stable $r+s$-pointed genus
$(g+1)$-curve by identifying $x_{r+s+1}$ with $O$. The normal bundle of this
stratum is just the exterior tensor product of the relative tangent bundles of
the factors. Denote the complement of its zero section by $U$ and let $U_e$ be
a general fibre of the projection of $U\to \M_1^1$. We can identify $U_e$ with
$\M _{g,1}^{r+s}$ and there is a natural restriction homomorphism
$$
H^{\bullet}(\M _{g+1}^{r+s};\Q )\to H^{\bullet}(U;\Q )\to H^{\bullet}(U_e;\Q
)\cong H^{\bullet}(\M _{g,1}^{r+s};\Q ).
$$
The $(\C ^*)^r$-bundles that lie over these spaces determine likewise a
homomorphism $H^{\bullet}(\M _{g+1,r}^s;\Q )\to H^{\bullet}(\M _{g,1+r}^s;\Q
)$. This is the one we were after.

Thus $H^{\bullet}(\G _{\infty };\Q )$ acquires a canonical mixed Hodge
structure. The composite map
$$
H^{\bullet}(\G _{\infty};\Q )\otimes A_X\to
H^{\bullet}(\G _g;\Q )\otimes A_X\cong H^{\bullet}(\M _g;\Q )\otimes A_X\to
H^{\bullet}(\Cu ^X _g;\Q )
$$
is evidently a morphism of mixed Hodge structures.

\medskip
For every $i\in X$ we have a projection $f_i:\Cu _g^X\to \Cu _g^{X-\{ i\}}$
and an obvious inclusion $A^{\bullet}_{X-\{i\}}\hookrightarrow A^{\bullet}_X$
(with image the linear
combinations of monomials in which the $u_I$ with $i\in I$ occur with exponent
$0$). The two are related:

\proclaim{\label Lemma} The map
$f_i^*:H^{\bullet}(\Cu _g^{X-\{ i\}};\Q ) \to H^{\bullet}(\Cu _g^X;\Q )$
is covered by the inclusion $A^{\bullet}_{X-\{i\}}\hookrightarrow
A^{\bullet}_X$
tensorized with the identity of $H^{\bullet}(\G _{\infty};\Q )$.
\endproclaim

The proof is straightforward.

Denote by $\overline{H}^{\bullet}(\Cu _g^X;\Q )$ the quotient of
$H^{\bullet}(\Cu _g^{X};\Q )$ by the subspace spanned by the images of $f_i^*,
u_if_i^*$, $i\in X$.

\proclaim{\label Corollary}
Put
$$
A'{}^{\bullet}_X:=\bigoplus _{P|X} (\prod _{\{ i\}\in P} u_i^2)\Q[u_I:I\in
P]a_P.
$$
Then there is a natural homomorphism of mixed Hodge structures
$$
H^{\bullet}(\G _{\infty};\Q )\otimes A'{}^{\bullet}_X\to\overline
{H}^{\bullet}(\Cu _g^X;\Q ) $$
which is an isomorphism in degree $\le N(g)$.
\endproclaim

\head
\section The Schur--Weyl functor
\endhead

\label We continue to denote by $X$ a finite nonempty set.
The product  $S_g^X$ comes with an obvious $\Sy _X$-action and so there is a
resulting action of $\Sy _X$ on $H^{\bullet}(S_g^X;\Q )$. Any  diffeomorphism
$f$
of $S_g$ induces a diffeomorphism of $S_g^X$ commuting with this action. Thus
is
obtained an action of the product $\Gamma _g\times \Sy _X$ on
$H^{\bullet}(S_g^X;\Q )$. As before, $V_g:=H^1(S_g;\Q )$ and $\symp
(V_g)$ denotes its symplectic group. It is easily seen that the $\Gamma
_g\times
\Sy _X$-action factorizes through one of $\symp (V_g)\times\Sy _X$.

Let us write the total
cohomology $H^{\bullet}(S_g;\Q )$ as $\Q \oplus V_gt\oplus \Q u$ where $u$ is
the
canonical class (we assume $g\ge 2$ here) and $t$ shifts the degree by $1$. The
K\"unneth rule gives an isomorphism
$$
H^{\bullet}(S_g^X;\Q )=\bigoplus _{I,J\subset X;I\cap J=\emptyset}
V_g^{\otimes I}t^Iu^J,
$$
where $t^I$ should be thought of as a generator of the signum representation of
$\Sy _I$ placed in degree $|I|$ and $u^J=\prod _{j\in J}u_j$. The
multiplication
in the right-hand side obeys the Koszul sign rule and is given by contractions
stemming from the symplectic form on $V_g$. If we define
$\overline{H}^{\bullet}(S_g^X;\Q )$ as in the relative case, that is, as
the quotient of $H^{\bullet}(S_g^X;\Q)$ by the span of the images of
$H^{\bullet}(S_g^{X-\{ i\}};\Q), u_iH^{\bullet}(S_g^{X-\{ i\}};\Q)$, $i\in X$,
then we see that in terms of the K\"unneth decomposition this reduces to
the single summand $V_g^{\otimes X}t^X$.

\medskip\label
Since $\pi :\Cu _g^X\to \M _g$ is a projective morphism
whose total space is an orbifold, it follows from a theorem of Deligne
\cite{1} that its Leray spectral sequence degenerates at the $E_2$-term over
$\Q$. In other words, $H^{\bullet}(\Cu _g^X;\Q )$ has a canonical
decreasing filtration $L^{\bullet}$, the {\it Leray filtration}, such that
there is a natural isomorphism
$$
\operatorname{Gr}_L^kH^n (\Cu _g^X;\Q )\cong H^k(\M _g;R^{n-k}\pi
_* \Q ).
$$
The maps $f_i^*$ in \refer{2.6} are strict with respect to
the Leray filtrations. So if we combine this with \refer{2.7} we find:

\proclaim{\label Corollary}
There is a natural graded $\Sy _s$-equivariant map of $H^{\bullet}(\G
_{\infty};\Q )$-modules
$$
H^{\bullet}(\G _{\infty};\Q )\otimes A'{}^{\bullet}_X
\to H^{\bullet}(\G _g; V_g ^{\otimes X})t^X
$$
which is an isomorphism in degree $\le N(g)$.
\endproclaim

We want to decompose $V_g^{\otimes X}$ as a $Sp(V_g)\times\Sy
_X$-representation. Following Weyl this is done in two steps.
The first step involves Weyl's representation $V_g^{\langle X\rangle}$ whose
definition we presently recall. Let $\omega\in V_g\otimes V_g$ correspond to
the
symplectic form on $V_g$. For every ordered pair $(i,j)\in X$ with $i\not= j$,
we have a natural homomorphism $V_g^{\otimes (X-\{i,j\})}\to V_g^{\otimes X}$
defined by placing $\omega$ in the $(i,j)$-slot. Reversing the order gives
minus
this map. So the natural assertion is that we have a map
$$
\bigoplus _{I\subset X; |I|=2} V_g^{\otimes (X -I)}t^I\to V_g^{\otimes X}.
$$
It is easy to see that this map is injective; its cokernel is by definition
$V_g^{\la X\ra}$. Notice that $V_g^{\langle
X\rangle}$ is in a natural way a representation of $\symp (V_g)\times\Sy _X$.
The second step is the decomposition of $V_g^{\langle X\rangle}$: Weyl proved
the remarkable fact that
$$
V_g^{\langle X\rangle}\cong \bigoplus _{\lambda} S_{\langle
{\lambda}\rangle}(V_g)\boxtimes (\lambda ),
$$
where $\lambda $ runs over the numerical partitions of $|X|$ in at most $g$
parts and $(\lambda )$ denotes the corresponding equivalence class of
irreducible representations of $\Sy _X$. In particular, all these irreducible
representations appear with multiplicity one. A modern account of the proof and
of related results can be found in the book by Fulton and Harris \cite{2}.

\proclaim{\label Theorem}
If $A''{}^{\bullet}_X\subset A'{}^{\bullet}_X$ is the $\Q [u_i:i\in
X]$-submodule defined by
$$
A''_X=\bigoplus _{P|X} (\prod _{\{ i\}\in I}u_i^2)(\prod _{I\in P;|I|=2}u_I)
\Q [u_I :I\in P]a_P,
$$
then there is a natural graded $\Sy _X$-equivariant map of $H^{\bullet}(\G
_g,\Q
)[u_i:i\in X]$-modules
$$
H^{\bullet}(\G _{\infty};\Q )\otimes A''{}^{\bullet}_X
\to H^{\bullet}(\G _g; V_g ^{\la X\ra })t^X,
$$
and this map is an isomorphism in degree $\le N(g)$.
\endproclaim

\demo{Proof} If $I\subset X$ is a two-element subset, then the Poincar\'e dual
of the hyperdiagonal $\Cu _g(P_I)\subset \Cu _g^X$ is $u_I$. The restriction of
this
element to the fiber $S_g^X$ is $(\sum _{i\in I}u_i)+\omega t^I$. So the map
$$
H^{\bullet -2}(\G _g; V_g ^{\otimes (X-I)})t^{X-I}\to
H^{\bullet}(\G _g; V_g ^{\otimes X})t^X
$$
defined by multiplication with $\omega t^I$ is covered by the map
$H^{\bullet}(\G _{\infty};\Q )\otimes A'{}^{\bullet}_{X-I}\to H^{\bullet}(\G
_{\infty};\Q )\otimes
A'{}^{\bullet}_X$ which is multiplication by $u_I-\sum_{i\in I}u_i$. The
theorem follows.
\enddemo

\medskip
\demo{Proof of \refer{1.1}} We decompose both members of the stable isomorphism
\refer{3.4} according to the action of $\Sy _X$. Let $\lambda$ be a numerical
partition of $|X|$ and let $\lambda '$ be the conjugate partition. According to
Weyl's decomposition theorem, the $(\lambda )$-isotypical component
of $H^{\bullet}(\G _g, V_g ^{\la X\ra})$ is $H^{\bullet}(\G _g,
S_{\la\lambda\ra}(V_g))$. Note that this is also the isotypical component of
type $\lambda '$ of $H^{\bullet}(\G _g; V_g ^{\langle X\rangle})t^X$. On the
other hand, it is easily seen that the $(\lambda ')$-isotypical component
of $A''{}^{\bullet}_X$ can be identified with $t^sB^{\bullet}(\lambda)$.

Since the Leray spectral sequence \refer{3.2} is a spectral sequence of mixed
Hodge structures, the homomorphism of \refer{1.1} is actually a morphism of
mixed Hodge structures. The theorem follows. \enddemo

\head
\section Stable cohomology of the universal Abel--Jacobi map
\endhead

If we give $S_g$ a complex structure, then $S_g$ becomes a compact Riemann
surface $C$ of genus $g$, so that we have defined an  Abel--Jacobi map $
\Sym ^s(C)\to\Pic ^s(C)$. The induced map on cohomology has been determined
by Macdonald:

\proclaim{\label Proposition}
{\rm (Macdonald \cite{7})} Identify $H^{\bullet}(\Pic ^s(C);\Q )$ with the
exterior algebra on $V_g$ so that the Abel--Jacobi map determines an algebra
homomorphism $\wedge ^{\bullet}V_g\to H^{\bullet}(\Sym ^s(C);\Q )$. Let $\wedge
^{\bullet}V_g[y]\to H^{\bullet}(\Sym ^s(C);\Q )$ be the extension that
sends the indeterminate $y$ of degree two to the sum of the fundamental classes
of the factors. Then this map is
surjective and its its kernel is the degree $>s$-part of the ideal $\Cal{I}$ in
$\wedge ^{\bullet}V_g [y]$ generated by $\{ v\wedge v'-(v.v')y : v,v'\in V_g
\}$.
\endproclaim

If $s>2g-2$, then the Abel--Jacobi map is the
projectivization of a vector bundle of rank $s+1-g$ over $\Pic (C)$ and we
can interpret the image of $y$ as the first Chern class of the associated line
bundle over $\Sym ^s(C)$. Macdonald also expresses the Poincar\'e
duals of the diagonals of $C^s$ in terms of this presentation.

\smallskip
It is our aim to make a corresponding discussion for the universal situation in
the stable range.

The morphism $\Cu _g\to \M _g$ defines
a relative Picard bundle $\Pic (\Cu _g/\M _g)\to \M_g$  (in the orbifold sense)
whose connected components are still indexed by the degree: $\Pic ^k(\Cu _g /
\M _g)$, $k\in\Z$. The degree $0$-component is called the {\it universal
Jacobian} and is also denoted $\J _g$.

\proclaim{\label Lemma}
For every $k\in\Z$, there is a natural isomorphism
$$
H^{\bullet}(\Pic ^k(\Cu _g/\M _g);\Q )\cong H^{\bullet}(\J _g;\Q ).
$$
\endproclaim
\demo{Proof} The relative canonical sheaf defines a section of $\Pic
^{2g-2}(\Cu
_g/\M _g)$. On a suitable Galois cover $\tilde{\M} _g\to \M _g$ (with Galois
group $G$, say) this section becomes divisible by $2g-2$ and thus produces a
section of  $\Pic ^1(\tilde\Cu _g/\tilde{\M } _g)$. This determines an
isomorphism $\Pic ^k(\tilde\Cu _g/\tilde\M _g)\cong \Pic ^0(\tilde\Cu
_g/\tilde\M _g)$. Although this isomorphism will not in general be
$G$-equivariant, it will differ from any $G$-translate by a section of $\Pic
^0(\tilde\Cu _g/\tilde M _g)$ of finite order, and so the induced map on
rational cohomology is $G$-equivariant. By passing to the $G$-invariants we
obtain an isomorphism
$H^{\bullet}(\Pic ^k(\Cu _g/\M _g);\Q )\cong H^{\bullet}(\J _g;\Q )$. One
checks
that this map does not depend on choices.
\enddemo

Let $X$ be a finite nonempty set as before. We wish to determine the subalgebra
of $\Sy _X$-invariants of $A^{\bullet}_X$, at least stably. Recall that an
additive basis of $A^{\bullet}_X$ consists of the set of elements of the form
$\prod _{I\in P} u_I^{r_I}$, where $P$ runs over the partitions of $X$ and
$r_I\ge |I|-1$. Let us define a partial ordering on the collection of
partitions
of $X$ by: $P\le Q$ if $P=Q$ or if for the smallest number $k$ such that the
$k$-element members of $P$ and $Q$ do not coincide every
$k$-element member of $P$ is a $k$-element member of $Q$. This determines a
partial ordering on the set of monomials: $\prod _{I\in P} u_I^{r_I}\le\prod
_{J\in Q} u_J^{s_J}$ if $P<Q$ or if $P=Q$ and $r_I\le s_I$ for all $I$. The
defining relations for $A^{\bullet}_X$ show that a product of two monomials
associated to partitions $P$ and $Q$ is a monomial associated to a partition
that dominates both $P$ and $Q$.

Denote by $S:=\sum _{\sigma\in\Sy _X}\sigma$ the symmetrizer operator (acting
on $A_X$).
Given a partition $P$ of $X$, then the smallest terms in
$\prod _{I\in P} S(u_I^{r_I})$ are monomials associated to a partition that is
a
$\Sy _X$-translate of $P$. In this expression the part associated to $P$
is the subsum corresponding to the partial symmetrizer $S_P:=\sum
_{\sigma\in\Sy
_X;\sigma (P)=P}\sigma$. If $l_k$ is the number of members of $P$ with $k$
elements, then the image of $S_P$ can be identified with
$$
\Q [c_1,c_2,\dots ,c_{l_1}]\otimes
\bigotimes _{k\ge 2;l_k>0} (c_{l_k}^{k-1}\Q [c_1,c_2,\dots ,c_{l_k}]).
$$
Here $c_l$ in the tensor factor with index $k$ is to be thought of as the $l$th
elementary symmetric function in the $u_I$'s with $I\in P$ and $|I|=k$ (and so
has degree $2l$); the appearance of $c_{l_k}^{k-1}$ comes from the condition
$r_I\ge |I|-1$.
So if we put
$$
C^{\bullet}_{\infty}:=
\Q [c_1,c_2,\cdots ]\otimes
\bigotimes _{k\ge 2}\bigl( \Q\oplus \bigoplus _{l\ge 1}c_l^{k-1}\Q
[c_1,c_2,c_3,\dots ,c_l]\bigr) ,
$$
then we find:

\proclaim{\label Corollary}
There is natural surjective homomorphism of graded algebra's
$C^{\bullet}_{\infty}\to (A^{\bullet}_X)^{\Sy _X}$.
Its restriction to the $\Q$-span of all the monomials involving the variables
$c^{(k_1)}_{l_1},\dots ,c^{(k_r)}_{l_r}$ (where $c_l^{(k)}$ denotes the
variable $c_l$ that occurs in the $k$th tensor power) with $\sum _i k_il_i\le
|X|$ is a linear isomorphism.
\endproclaim

If $Y$ is another finite set with $|Y|\ge |X|$, then any injection
$f: X\hookrightarrow Y$ induces an algebra homomorphism $A^{\bullet}_X\to
A^{\bullet}_Y$. Since all such injections are in the same $\Sy _Y$-orbit
they give rise to the same algebra homomorphism $(A^{\bullet}_X)^{\Sy
_X}\to (A^{\bullet}_Y)^{\Sy _Y}$. In particular,
$(\tilde A^{\bullet}_X)^{\Sy _X}$ only depends on $|X|$. So if we write
$\tilde C^{\bullet}_{|X|}$ for this algebra, then we have a direct system
$\cdots\to C^{\bullet}_s\to C^{\bullet}_{s+1}\to\cdots $. The corollary shows
that the limit of this direct system can be identified with
$C^{\bullet}_{\infty}$.
It follows from \refer{2.3} that we have an algebra homomorphism
$$
H^{\bullet}(\G _{\infty};\Q )\otimes C^{\bullet}_s\to
H^{\bullet}((\Cu _g ^s);\Q )^{\Sy _s}\cong H^{\bullet}((\Cu
_g ^s)^{\Sy _s};\Q )
$$
that is an isomorphism in degree $\le N(g)$. In the limit this yields a
homomorphism
$$
H^{\bullet}(\G _{\infty};\Q )\otimes C^{\bullet}_{\infty}\to H^{\bullet}((\Cu
_g ^s)^{\Sy _s};\Q )
$$
that is an isomorphism in degree $\le \min (2s,N(g))$. The image of
$c_1\otimes 1\otimes 1\otimes\cdots$ is easily seen to be proportional to the
element $y$ appearing in Macdonald's theorem. We put
$$
C'{}^{\bullet}_{\infty}:=
\Q [c_2,c_3,\cdots ]\otimes\bigotimes _{k\ge 2}\bigl( \Q+\oplus _{l\ge
1}c_l^{k-1}\Q [c_1,c_2,c_3,\dots ,c_l]\bigr) .
$$

\proclaim{\label Theorem} The algebra homomorphism above fits in a commutative
square of algebra homomorphisms
$$ \matrix
H^{\bullet}(\G _{\infty};\Q )\otimes C^{\bullet}_{\infty}&\longrightarrow &
H^{\bullet}((\Cu _g ^s)^{\Sy _s};\Q )\\
\cup &&\uparrow\\
H^{\bullet}(\G _{\infty};\Q )\otimes C'{}^{\bullet}_{\infty}&\longrightarrow &
H^{\bullet}(\Pic ^s(\Cu _g/\M _g);\Q ),
\endmatrix
$$
in which the right vertical map is induced by the Abel--Jacobi map. The lower
horizontal map is an isomorphism in degree $\le \min (s,N(g))$ so that in the
limit we have an isomorphism
$$
H^{\bullet}(\G _{\infty};\Q )\otimes C'{}^{\bullet}_{\infty}\cong
\bigoplus _{s=0}^{\infty} H^{\bullet}(\G _{\infty};\wedge ^s)t^s.
$$
\endproclaim
\demo{Proof}
If we combine Macdonald's theorem with the Leray spectral of sequence of the
map $(\Cu _g ^s)^{\Sy _s}\to \M_g$, then we see that the map
$$
H^{\bullet}(\Pic ^s(\Cu _g/\M _g);\Q )[y]\to
H^{\bullet}((\Cu _g ^s)^{\Sy _s};\Q )
$$
that sends $y$ to $c_1\otimes 1\otimes 1\otimes\cdots $ is an isomorphism in
degree $\le s$. The theorem follows from this.
\enddemo

\enddocument
\bye